\documentclass[11pt,a4paper]{article}

\usepackage{url}

\usepackage{latexsym,amssymb, amsthm}
\usepackage[centertags]{amsmath}
\usepackage{epsfig,pifont}
\usepackage[inactive]{srcltx}  

\usepackage{graphics,helvet,times,mathptm,epic,eepic}

\usepackage{color}

\usepackage{newlfont}

\usepackage{amsfonts,bbm}

\def\beq{\begin{equation}}
\def\eeq{\end{equation}}

\theoremstyle{remark}
\newtheorem{example}{Example}[section]

\def\G1{\hbox{$\displaystyle{\mbox{\ding{172}}}$}}

\begin{document}

%

%
\title{Methodology of Numerical Computations\\ with
Infinities and Infinitesimals}

\newcommand{\nms}{\normalsize}
\author{  {   \bf Yaroslav D. Sergeyev\footnote{Yaroslav D.
Sergeyev, Ph.D., D.Sc., is Distinguished Professor at the University
of Calabria, Rende, Italy.
 He is also Full Professor (a part-time contract) at the N.I.~Lobatchevsky State University,
  Nizhni Novgorod, Russia and Affiliated Researcher at the Institute of High Performance
  Computing and Networking of the National Research Council of Italy.   }
  \footnote{This paper presents Honorable Lagrange Lecture delivered in Turin on 16 April 2010. This research was partially supported by the
Russian Federal Program ``Scientists and Educators in Russia of
Innovations", contract number 02.740.11.5018.}
    }\\ \\ [-2pt]
      \nms Dipartimento di Elettronica, Informatica e Sistemistica,\\[-4pt]
       \nms   Universit\`a della Calabria,\\[-4pt]
       \nms 87030 Rende (CS)  -- Italy\\ \\[-4pt]
       \nms http://wwwinfo.deis.unical.it/$\sim$yaro\\[-4pt]
         \nms {\tt  yaro@si.deis.unical.it }
}

\date{}

\maketitle

%

\begin{abstract}
A recently developed computational  methodology for executing
numerical calculations
 with infinities and infinitesimals is described in
this paper. The developed approach has a pronounced applied
character and is based on the principle `The part is less than the
whole' introduced by Ancient Greeks. This principle is used with
respect to all numbers (finite, infinite, and infinitesimal) and
to all sets and processes (finite and infinite). The point of view
on infinities and infinitesimals (and in general, on Mathematics)
presented in this paper uses strongly physical  ideas emphasizing
interrelations holding between a
  mathematical object under the observation and   tools used for
this observation. It is shown   how a new numeral system allowing
one to express different infinite and infinitesimal quantities  in
a unique framework can be used for theoretical and computational
purposes. Numerous examples dealing with infinite sets, divergent
series, limits, and probability theory are given.
\end{abstract}

{\bf AMS Subject Classification:} 03E65, 65-02, 65B10, 60A10 \\

\section{Introduction}
\label{s1}

The concept of infinity attracted the attention of people during
millenniums  (see monographs
\cite{Cantor,Conway,Godel,Hardy,Hilbert,Leibniz,Mayberry,
Newton,Robinson} and references given therein). To emphasize
importance of the subject for modern Mathematics it is sufficient
to mention that the Continuum Hypothesis related to infinity has
been included by David Hilbert as the Problem Number One in his
famous list of 23 unsolved mathematical problems (see
\cite{Hilbert}) that have influenced strongly development of
Mathematics in the XX-th century.

There exist different ways to generalize traditional arithmetic
for finite numbers to the case of infinities and infinitesimals
(see, e.g., \cite{Cantor,Conway,Robinson}  and references given
therein). However, the arithmetics developed for infinite numbers
up to now were quite different with respect to the finite
arithmetic we are used to deal with. Very often they leave
undetermined many operations where infinity takes part (for
example, $\infty-\infty$, $\frac{\infty}{\infty}$, sum of
infinitely many items, etc.) or use a representation of infinite
numbers based on infinite sequences of finite numbers. In spite of
these crucial difficulties and due to enormous importance of the
concept of infinity in science, people try to introduce infinity
in their work with computers. We can mention the IEEE Standard for
Binary Floating-Point Arithmetic containing representations for
$+\infty$ and $-\infty$ and incorporation of these notions in the
interval analysis implementations.

The development of the modern views on infinity and infinitesimals
strangely enough was not simultaneous. The point of view on
infinity accepted nowadays takes its origins from the famous ideas
of Georg Cantor (see \cite{Cantor}) who has shown that there exist
infinite sets having different cardinalities. On the other hand,
in the early history of Calculus, arguments involving
infinitesimals played a pivotal role in the differential Calculus
developed by Leibniz and Newton (see \cite{Leibniz,Newton}). The
notion of an infinitesimal, however, lacked a precise mathematical
definition and, in order to provide a more rigorous foundation for
the Calculus, infinitesimals were gradually replaced by the
d'Alembert--Cauchy concept of a limit (see
\cite{Cauchy,DAlembert}).

The creation of a rigorous mathematical theory of infinitesimals
  remained an
open problem until the end of the 1950s when Robinson (see
\cite{Robinson}) has introduced his famous non-standard Analysis
approach. He has shown that non-archimedean ordered field
extensions of the reals contained numbers that could serve the
role of infinitesimals and their reciprocals   could serve as
infinitely large numbers. Robinson then has derived the theory of
limits,  and more generally of Calculus, and has found a number of
important applications of his ideas in many other fields of
Mathematics (see \cite{Robinson}).

In his approach,  Robinson used mathematical tools and terminology
(cardinal numbers, countable sets, continuum, one-to-one
correspondence, etc.) taking their origins from the  ideas of
Cantor (see \cite{Cantor}) introducing so all advantages and
disadvantages of Cantor's theory   in non-standard Analysis, as
well. In fact, it is well known nowadays that while dealing with
infinite sets, Cantor's approach leads to some counterintuitive
situations that often are called by non-mathematicians as
`paradoxes'. For example, the set of even numbers, $\mathbb{E}$,
can be put in a one-to-one correspondence with the set of all
natural numbers, $\mathbb{N}$, in spite of the fact that
$\mathbb{E}$ is  a proper subset of $\mathbb{N}$:
 \beq
\begin{array}{lccccccc}
  \mbox{even numbers:}   & \hspace{5mm} 2, & 4, & 6, & 8,  & 10, & 12, & \ldots    \\

& \hspace{5mm} \updownarrow &  \updownarrow & \updownarrow  & \updownarrow  & \updownarrow  &  \updownarrow &   \\

  \mbox{natural numbers:}& \hspace{5mm}1, &  2, & 3, & 4 & 5,
       & 6,  &    \ldots \\
     \end{array}
\label{4.4.1}
 \eeq
In contrast, we can observe that for finite sets, if a set $A$ is
a proper subset of a set $B$ then it follows that the number of
elements of the set $A$ is smaller than the number of elements of
the set $B$.

Another famous example that is  difficult for understanding  for
many people   is Hilbert's paradox of the Grand Hotel having the
following formulation. In a normal hotel with a finite number of
rooms no more new guests can be accommodated if it is full.
Hilbert's Grand Hotel has an infinite number of rooms (of course,
the number of rooms is countable, because the rooms in the Hotel
are numbered). Due to Cantor, if a new guest arrives at the Hotel
where every room is occupied, it is, nevertheless, possible to
find a room for him. To do so, it is necessary to move the guest
occupying room 1 to room 2, the guest occupying room 2 to room 3,
etc. In such a way room 1 will be ready for the newcomer and, in
spite of our assumption that there are no available rooms in the
Hotel, we have found one.

These results are very difficult to be fully realized by anyone
who is not a mathematician since in our every day experience in
the world around us the part is always less than the whole and if
a hotel is complete there are no places in it. In order to
understand how it is possible to tackle the situations  discussed
above   in accordance with the principle `the part is less than
the whole' let us consider a study published in \textit{Science}
(see \cite{Gordon}) where the author describes a primitive tribe
living in Amazonia - Pirah\~{a} - that uses a very simple numeral
system\footnote{ We remind that \textit{numeral} is a symbol or
group of symbols that represents a \textit{number}. The difference
between numerals and numbers is the same as the difference between
words and the things they refer to. A \textit{number} is a concept
that a \textit{numeral} expresses. The same number can be
represented by different numerals. For example, the symbols `10',
`ten', and `X' are different numerals, but they all represent the
same number.} for counting: one, two, many. For Pirah\~{a}, all
quantities larger than two are just `many' and such operations as
2+2 and 2+1 give the same result, i.e., `many'. Using their weak
numeral system Pirah\~{a} are not able to see, for instance,
numbers 3, 4, 5, and 6, to execute arithmetical operations with
them, and, in general, to say anything about these numbers because
in their language there are neither words nor concepts for that.
Moreover, the weakness of Pirah\~{a}'s numeral system leads to
such results as
\[
\mbox{`many'}+ 1= \mbox{`many'},   \hspace{1cm}    \mbox{`many'} +
2 = \mbox{`many'},
\]
which are very familiar to us  in the context of views on infinity
used in the traditional calculus
\[
\infty + 1= \infty,    \hspace{1cm}    \infty + 2 = \infty.
\]
Thus, the modern mathematical numeral systems allow us to
distinguish a larger quantity of finite numbers with respect to
Pirah\~{a} but give  similar results  when we speak  about
infinite numbers.

The arithmetic of Pirah\~{a} involving the numeral `many' has also
a clear similarity with the arithmetic proposed by Cantor for his
Alephs. This similarity becomes even stronger if one considers
another Amazonian tribe -- Munduruk\'u (see \cite{Pica}) -- who
fail in exact arithmetic with numbers larger than 5 but are able
to compare and add large approximate numbers that are far beyond
their naming range. Particularly, they use the words `some, not
many' and `many, really many' to distinguish two types of large
numbers (in this connection think about Cantor's $\aleph_0$ and
$\aleph_1$).

These observations lead  us to the following idea:
\textit{Probably our difficulty in working with infinity is not
connected to the nature of infinity but is a result of inadequate
numeral systems used to express infinite numbers.} Analogously,
Pirah\~{a} are not able to distinguish numbers 3 and 4 not due to
the nature of these numbers but due to the weakness of the numeral
system that Pirah\~{a} use.

In this paper, we show how the introduction of a new numeral
 allows one to express different infinite and
infinitesimal quantities. Taken together with a new (physically
oriented) methodology for Mathematics, the new numeral system
 can be used for theoretical and computational purposes
 using the Infinity Computer (see
\cite{Sergeyev_patent}) able to work numerically with infinite and
infinitesimal numbers expressed in the new numeral system.

\section{From absolute mathematical truths to  their  relativity and accuracy of mathematical results}
\label{s2}

 In this section, we give   a brief
introduction to the new methodology    that can be found in a
rather comprehensive form in the survey \cite{informatica}
downloadable from \cite{www} (see also the monograph
\cite{Sergeyev} written in a popular manner and
\cite{Dif_Calculus} describing the foundations of a new
differential calculus). Numerous examples of the usage of the
proposed methodology can be found in \cite{Sergeyev,chaos,
Menger,Korea,first,Sergeyev_Garro,www}. The goal of the entire
operation is to propose a way of thinking that would allow us to
work with finite, infinite, and infinitesimal numbers in the same
way, namely, in the way we are used to deal with finite quantities
in the world around us.

In order to start, let us make some observations. As   was
mentioned above, foundations of the modern Set Theory dealing with
infinity have been developed starting from the end of the XIX-th
century until more or less the first decades of the XX-th century.
Foundations of the classical Analysis dealing both with infinity
and infinitesimals have been developed even earlier, more than 200
years ago. The goal of its creation was to produce mathematical
tools allowing one to solve problems arising in the real world in
that time. As a result, classical Analysis was build using the
common in that time background of ideas that people had about
Physics (and Philosophy). Thus, these parts of Mathematics do not
include numerous achievements of Physics of the XX-th century. In
fact, the classical Analysis operates with  absolute truths   and
ideas of relativity and quants are not reflected in it. Let us
give just one example to clarify this point.

In modern Physics, the `continuity'  of an object is relative. If
we observe a table by eye, then we see it as being continuous. If
we use a microscope for our observation, we see that the table is
discrete. This means that \textit{we decide} how to see the
object, as a continuous or as a discrete, by the choice of the
instrument for the observation. A weak instrument -- our eyes --
is not able to distinguish its internal small separate parts
(e.g., molecules) and we see the table as a continuous object. A
sufficiently strong microscope allows us to see the separate parts
and the table becomes discrete but each small part now is viewed
as continuous.

In contrast, in   traditional Mathematics, any mathematical object
is either continuous or discrete. For example, the same function
cannot be   both  continuous and discrete. Thus, this
contraposition of discrete and continuous in the traditional
Mathematics does not reflect  properly the physical situation that
we observe in practice.

Note that even results of Robinson    made in the middle of the
XX-th century   have been also directed to a reformulation of the
classical Analysis   in terms of infinitesimals and not to the
creation of a new kind of Analysis that would incorporate new
achievements of Physics. In fact, he wrote in paragraph 1.1 of his
famous book \cite{Robinson}: `It is shown in this book that
Leibniz' ideas can be fully vindicated and that they lead to a
novel and fruitful approach to classical Analysis and to many
other branches of mathematics'.

In order to overcome this delay with the introduction of ideas of
Physics of the XX-th century in Mathematics, the point of view on
infinities and infinitesimals (and in general, on Mathematics)
presented in this paper uses strongly relativity and
interrelations holding between the object of an observation and
the tool used for this observation. The latter is directly related
to connections between numeral systems used to describe
mathematical objects and the objects themselves. Numerals that we
use to write down numbers, functions, etc.  are among our tools of
the investigation and, as a result, they strongly influence our
capabilities to study mathematical objects.

This separation (having an evident physical spirit) of
mathematical objects from tools used for their description is
crucial for our study but it is used rarely in contemporary
Mathematics. In fact, the idea of finding an adequate (absolutely
the best) set of axioms for one or another field of Mathematics
continues to be among the most attractive goals for contemporary
mathematicians. Usually, when it is necessary to define a concept
or an object, logicians try to introduce a number of axioms
\textit{defining} the object. However, this way is fraught with
danger because of the following reasons.

First, when we describe a mathematical object or concept we are
limited by the expressive capacity of the language we use to make
this description. A   richer language allows us to say more about
the object and a weaker language -- less. Thus, development of the
mathematical (and not only mathematical) languages leads to a
continuous necessity of a transcription and specification of
axiomatic systems. Second, there is no   guarantee that the chosen
axiomatic system defines `sufficiently well' the required concept
and a continuous comparison with practice is required in order to
check the goodness of the accepted set of axioms. However, there
cannot be again any guarantee that the new version will be the
last and definitive one. Finally, the third limitation already
mentioned above    has been discovered by G\"odel in his two
famous incompleteness theorems (see~\cite{Godel_1931}).

It should be emphasized that in Linguistics, the relativity of the
language with respect to the world around is a well known thing.
It has been formulated in the form of the Sapir--Whorf thesis
(see~\cite{Whorf,Sapir}) also known as the `linguistic relativity
thesis'. As becomes  clear from  its name, the thesis does not
accept the idea of the universality of language and postulates
that the nature of a particular language influences the thought of
its speakers. The thesis challenges the possibility of perfectly
representing the world with language, because it implies that the
mechanisms of any language condition the thoughts of its speakers.

Thus, our point of view on axiomatic systems is  different. It is
significantly more applied and less ambitious and is related only
to utilitarian necessities to make calculations.    In contrast to
the modern mathematical fashion that tries to make all axiomatic
systems more and more precise (decreasing so degrees of freedom of
the studied part of Mathematics), we just define a set of general
rules describing how practical computations should be executed
leaving so as much space as possible for further, dictated by
practice, changes and developments of the introduced mathematical
language. Speaking metaphorically, we prefer to make a hammer and
to use it instead of describing what is a hammer and how it works.

Since our point of view on the mathematical world is significantly
more physical and more applied than the traditional one, it
becomes necessary to clarify it better. Let us formulate three
methodological postulates that will guide our further study and
will show where our positions are different with respect to the
tradition.

Traditionally, when mathematicians deal with infinite objects
(sets or processes) it is supposed   that human beings are able to
execute certain operations infinitely many times (e.g., see
(\ref{4.4.1})). However, since we live in a finite world and all
human beings and/or computers are forced to finish operations that
they have started, this supposition is not adopted.

 \textbf{Postulate 1.} \textit{There  exist
infinite and infinitesimal objects but   human beings and machines
are able to execute only a finite number of operations.}

Due to this Postulate, we accept a priori that we shall never be
able to give a complete description of infinite processes and sets
due to our finite capabilities.

The second postulate is adopted  following the way of reasoning
used in natural sciences where researchers use tools to describe
the object of their study and the used instrument   influences the
results of the observations. When a physicist uses a weak lens $A$
and sees two black dots in his/her microscope he/she does not say:
The object of the observation \textit{is} two black dots. The
physicist is obliged to say: the lens used in the microscope
allows us to see two black dots and it is not possible to say
anything more about the nature of the object of the observation
until we  change the instrument - the lens or the microscope
itself - by a more precise one. Suppose that  he/she changes the
lens and uses a stronger lens $B$ and is able to observe that the
object of the observation is viewed as ten (smaller) black dots.
Thus, we have two different answers: (i) the object is viewed as
two dots if the lens $A$ is used; (ii) the object is viewed as ten
dots by applying the lens $B$. Which of the answers is correct?
Both. Both answers are correct but with the different accuracies
that depend on the lens used for the observation.

The same happens in Mathematics studying natural phenomena,
numbers, and objects that can be constructed by using numbers.
Numeral systems used to express numbers are among the instruments
of observations used by mathematicians. The usage of powerful
numeral systems gives the possibility to obtain more precise
results in Mathematics in the same way as usage of a good
microscope gives the possibility of obtaining more precise results
in Physics. However, the capabilities of the tools will be always
limited due to Postulate 1 (we are able to write down only a
finite number of symbols when we wish to describe a mathematical
object) and due to Postulate~2 we shall never tell, \textbf{what
is}, for example, a number but shall just observe it through
numerals expressible in a chosen numeral system.

 \textbf{Postulate
2.} \textit{We shall not   tell \textbf{what are} the mathematical
objects we deal with; we just shall construct more powerful tools
that will allow us to improve our capacities to observe and to
describe properties of mathematical objects.}

This means that mathematical results are not absolute, they depend
on mathematical languages used to formulate them, i.e., there
always exists an accuracy of the description of a mathematical
result, fact, object, etc. For instance, the result of Pirah\~{a}
$2+2=$ `many' is not wrong, it is just inaccurate. The
introduction of a stronger tool (in this case, a numeral system
that contains a numeral for a representation of the number four)
allows us to have a more precise answer.

It is necessary to comment upon another important aspect of the
distinction between a mathematical object and a mathematical tool
used to observe this object. The Postulates 1 and 2 impose us to
think always about \textit{the possibility to execute} a
mathematical operation by applying a numeral system. They tell us
that there always exist situations where we are not able to
express the result of an operation. Let us consider, for example,
the operation of construction of the successive element widely
used in number and set theories. In the traditional Mathematics,
the aspect  whether this operation can be executed is not taken
into consideration, it is supposed that it is always possible to
execute the operation $k=n+1$ starting from any integer $n$. Thus,
there is no any distinction between the existence of the number
$k$ and the possibility to execute the operation $n+1$ and to
express its result, i.e. to have a numeral that can express $k$.

Postulates 1 and 2 emphasize this distinction and tell us that:
(i)~ in order to execute the operation it is necessary to have a
numeral system allowing one to express both numbers, $n$ and $k$;
(ii)~for any numeral system there always exists a number $k$ that
cannot be expressed in it. For instance, for Pirah\~{a} $k=3$, for
Munduruk\'u $k=6$. Even for modern powerful numeral systems there
exist such a number $k$ (for instance, nobody is able to write
down a numeral in the decimal positional system having $10^{100}$
digits). Hereinafter we shall always emphasize the triad --
researcher, object of the investigation, and tools used to observe
the object -- in various mathematical and computational contexts
paying a special attention to the accuracy of the obtained
results.

Particularly, Postulate~2 means that, from our point of view,
axiomatic systems \textit{do not define} mathematical objects but
just determine formal rules for operating with certain numerals
reflecting some properties of the studied mathematical objects
using a certain mathematical language $L$. We are aware that the
chosen language $L$ has its accuracy and there always can exist a
richer language $\tilde{L}$ that would allow us to describe the
studied object better. Due to Postulate~1, any language has a
limited expressibility, in particular, there always exist
situations where the accuracy of the answers  expressible in this
language is not sufficient. Such situations lead to `paradoxes'
showing the boundaries of the applicability of a language (theory,
concept, etc.)

Let us return again to Pirah\~{a} and illustrate this point by
using their answers $2+1=$ `many' and $2+2=$ `many'. From these
two identities one can obtain the result $2+1= 2+2$ being a
`paradox'. From our point of view, this situation just determines
the boundaries of the applicability of their numeral system.

Finally, we adopt the principle of  Greeks mentioned above as  the
third postulate.

\textbf{Postulate 3.} \textit{The principle `The part is less than
the whole' is applied to all numbers (finite, infinite, and
infinitesimal) and to all sets and processes (finite and
infinite).}

Due to this declared applied statement, it becomes clear that the
subject of this paper is out of Cantor's approach and, as a
consequence, out of non-standard Analysis of Robinson. Such
concepts as bijection, numerable and continuum sets, cardinal and
ordinal numbers cannot be used in this paper because they belong
to the theory working with different assumptions. However, the
approach used here does not contradict Cantor and Robinson. It can
be viewed just as a more strong lens of a mathematical microscope
that allows one to distinguish more objects and to work with them.

\section{An infinite unit of measure expressible by a new numeral} \label{s3}

In \cite{Sergeyev,informatica}, a
  new  numeral system has
been developed in accordance with methodological Postulates~1--3.
It gives a possibility to execute numerical computations not only
with finite numbers but also with infinite and infinitesimal ones.
The main idea consists of the possibility to measure infinite and
infinitesimal quantities by different (infinite, finite, and
infinitesimal) units of measure.

A new infinite unit of measure   has been introduced for this
purpose   as the number of elements of the set $\mathbb{N}$ of
natural numbers. The new number is  called \textit{grossone} and
is expressed by the numeral \G1. It is necessary to stress
immediately that \G1 is neither Cantor's $\aleph_0$ nor $\omega$.
Particularly, it has both cardinal and ordinal properties as usual
finite natural numbers (see \cite{informatica}). Note also that
since \G1, on the one hand, and $\aleph_0$  (and $\omega$), on the
other hand, belong to different mathematical languages working
with different theoretical assumptions, they cannot be used
together. Analogously, it is not possible to use together Piraha's
`many' and the modern numeral 4.

Formally, grossone is introduced as a new number by describing its
properties postulated by the \textit{Infinite Unit Axiom} (IUA)
(see \cite{Sergeyev,informatica}). This axiom is added to axioms
for real numbers similarly to addition of the axiom determining
zero to axioms of natural numbers when integer numbers are
introduced. It is important to emphasize that we speak about
axioms for real numbers in sense of Postulate~2, i.e., axioms do
not define real numbers, they just define formal rules of
operations with  numerals in given  numeral systems (tools of the
observation)  reflecting so certain (not all) properties of the
object of the observation, i.e., properties of real numbers.

Inasmuch as it has been postulated that grossone is a number,  all
other axioms for numbers hold for it, too. Particularly,
associative and commutative properties of multiplication and
addition, distributive property of multiplication over addition,
existence of   inverse  elements with respect to addition and
multiplication hold for grossone as for finite numbers. This means
that  the following relations hold for grossone, as for any other
number
 \beq
 0 \cdot \G1 =
\G1 \cdot 0 = 0, \hspace{3mm} \G1-\G1= 0,\hspace{3mm}
\frac{\G1}{\G1}=1, \hspace{3mm} \G1^0=1, \hspace{3mm}
1^{\mbox{\tiny{\G1}}}=1, \hspace{3mm} 0^{\mbox{\tiny{\G1}}}=0.
 \label{3.2.1}
       \eeq
The introduction of the new numeral allows us to use it for
construction of various new numerals expressing infinite and
infinitesimal numbers  and to operate with them as with usual
finite constants. As a consequence, the numeral $\infty$   is
excluded from our new mathematical language (together with
numerals $\aleph_0, \aleph_1, \ldots$  and $\omega$). In fact,
since we are able now to express explicitly different infinite
numbers, records of the type $\sum_{i=1}^{\infty}a_i$ become  a
kind of $\sum_{i=1}^{many}a_i$, i.e., they are not sufficiently
precise. It becomes necessary not only to say that $i$ goes to
infinity, it is necessary to indicate to which point in infinity
(e.g., $\G1, 5\G1-1, \G1^2+3$, etc.) we want to sum up. Note that
for sums having a finite number of items the situation is the
same: it is not sufficient to say that the number of items in the
sum is finite, it is necessary to indicate explicitly the number
of items in the sum.

The appearance of new numerals expressing infinite and
infinitesimal numbers   gives us a lot of new possibilities. For
example, it becomes possible to develop a Differential Calculus
(see \cite{Dif_Calculus}) for functions that can assume finite,
infinite, and infinitesimal values and can be defined over finite,
infinite, and infinitesimal domains avoiding indeterminate forms
and divergences (all these concepts just do not appear in the new
Calculus). This approach allows us to work with derivatives and
integrals that can assume not only finite but infinite and
infinitesimal values, as well. Infinite and infinitesimal numbers
are not auxiliary entities in the new Calculus, they are full
members in it and can be used in the same way as finite constants.

Let us comment upon the nature of grossone and give some examples
illustrating its usage and, in particular, its direct links with
infinite sets.

\begin{example}
\label{e1_Lagrange}  Grossone has been introduced as  the number
of elements of the set $\mathbb{N}$ of natural numbers. As a
consequence, similarly  to the set
 \beq
  A=\{1, 2, 3, 4, 5\}
\label{4.1.deriva_0}
 \eeq
   consisting of
5 natural numbers where 5 is the largest number in $A$, \G1 is the
largest    number\footnote{This fact is one of the important
methodological differences with respect to non-standard analysis
theories where it is supposed that infinite numbers   do not
belong to $\mathbb{N}$.} in $\mathbb{N}$ and $\G1 \in \mathbb{N}$
analogously to the fact that 5 belongs to $A$. Thus, the set,
$\mathbb{N}$, of natural numbers can be written  in the form
 \beq
\mathbb{N} = \{ 1,2,  \hspace{3mm} \ldots  \hspace{3mm}
\frac{\G1}{2}-2, \frac{\G1}{2}-1, \frac{\G1}{2}, \frac{\G1}{2}+1,
\frac{\G1}{2}+2, \hspace{3mm}  \ldots \hspace{3mm} \G1-2,
\hspace{2mm}\G1-1, \hspace{2mm} \G1 \}. \label{4.1}
       \eeq
Note that traditional numeral systems did not allow us to see
infinite natural numbers
 \beq \ldots  \hspace{3mm}
\frac{\G1}{2}-2, \frac{\G1}{2}-1, \frac{\G1}{2}, \frac{\G1}{2}+1,
\frac{\G1}{2}+2, \hspace{3mm} \ldots  \hspace{3mm} \G1-2, \G1-1,
\G1. \label{4.1.deriva_1}
 \eeq
Similarly,
  Pirah\~{a}     are not able to see  finite numbers larger than 2
using their weak numeral system but these numbers are visible if
one uses a more powerful numeral system. Due to Postulate~2, the
same object  of observation -- the set $\mathbb{N}$ --   can be
observed by different instruments -- numeral systems -- with
different accuracies allowing one to express  more or less natural
numbers. \hfill $\Box$
\end{example}

This example illustrates also the fact that when we speak about
sets (finite or infinite) it is necessary to take care about tools
used to describe a set (remember Postulate~2). In order to
introduce a set, it is necessary to have a language (e.g., a
numeral system) allowing us to describe its elements and to
express the number of the elements in the set. For instance, the
set $A$ from (\ref{4.1.deriva_0}) cannot be defined using the
mathematical language of Pirah\~{a}.

Analogously, the words `the set of all finite numbers' do not
define a set completely from our point of view, as well. It is
always necessary to specify which instruments are used to describe
(and to observe) the required set and, as a consequence, to speak
about `the set of all finite numbers expressible in a fixed
numeral system'. For instance, for Pirah\~{a} `the set of all
finite numbers'  is the set $\{1, 2 \}$ and for Munduruk\'u  `the
set of all finite numbers' is the set $A$ from
(\ref{4.1.deriva_0}). As it happens in Physics, the instrument
used for an observation bounds the possibility of the observation.
It is not possible to say how we shall see the object of our
observation if we have not clarified which instruments will be
used to execute the observation.

\begin{example}
\label{e2_Lagrange} Infinite numerals constructed using \G1  allow
us to observe various infinite integers being   the number of
elements of infinite sets. For example, $\G1-1$ is the number of
elements of a set $B=\mathbb{N}\backslash\{b\}$, $b \in
\mathbb{N}$, and $\G1+1$ is the number of elements of a set
$A=\mathbb{N}\cup\{a\}$, where $a \notin \mathbb{N}$.

Due to Postulate~3, positive integers that are larger than
grossone do not belong to $\mathbb{N}$. However, numerals
expressing such numbers can be easily constructed and it can be
shown that they represent the number of elements of certain
infinite sets. For instance, $\G1^2$ is the number of elements of
the set
  $V$ of couples of natural numbers
  $$ V  =
\{ (a_1, a_2)  : a_1 \in   \mathbb{N}, a_2 \in   \mathbb{N}  \}.
 $$
By increasing $a_1$ and $a_2$ from 1 to \G1 we are able to write
down initial and final couples forming this
 set:
\[
  \begin{array}{ccccc}
(1,1), & (1,2), & \ldots & (1,\G1-1), & (1,\G1), \\
  (2,1), & (2,2), & \ldots & (2,\G1-1), & (2,\G1), \\
\ldots & \ldots & \ldots & \ldots & \ldots \\
\hspace{1mm}(\G1-1,1),  \hspace{1mm} & \hspace{1mm}  (\G1-1,2),\hspace{1mm} & \hspace{1mm} \ldots \hspace{1mm} & \hspace{1mm}(\G1-1,\G1-1), \hspace{1mm} & (\G1-1,\G1), \\
(\G1,1), & (\G1,2), & \ldots & (\G1,\G1-1), & (\G1,\G1).
  \end{array}
\]
 Analogously, the number $2^{\tiny{\G1}}$ is the number
of elements of the set
$$ U  =
\{ (a_1, a_2, \ldots a_{\tiny{\G1}-1}, a_{\tiny{\G1}})  : a_1 \in
\{1,2\}, a_2  \in \{1,2\}, \ldots a_{\tiny{\G1}-1} \in \{1,2\},
a_{\tiny{\G1}} \in \{1,2\} \}     $$
  and the number $\G1^{\tiny{\G1}}$ is the number
of elements of the set
$$ \hspace{15mm} W  =
\{ (a_1, a_2, \ldots a_{\tiny{\G1}-1}, a_{\tiny{\G1}})  : a_1 \in
\mathbb{N}, a_2  \in \mathbb{N},  \ldots a_{\tiny{\G1}-1} \in
\mathbb{N}, a_{\tiny{\G1}} \in \mathbb{N} \}. \hspace{10mm}  \Box
$$

\end{example}

As was mentioned above, the introduction of  grossone gives us a
possibility to compose new (in comparison with traditional numeral
systems) numerals and to see through them not only numbers
(\ref{4.1.deriva_0}) but also certain numbers larger than \G1. We
can speak about the set of \textit{extended natural numbers}
(including $\mathbb{N}$ as a proper subset) indicated as
$\widehat{\mathbb{N}}$ where
 \beq
  \widehat{\mathbb{N}} = \{
1,2, \ldots ,\G1-1, \G1, \G1+1, \G1+2, \G1+3, \ldots , \G1^2-1,
\G1^2. \G1^2+1, \ldots \}. \label{4.2.2}
       \eeq
The number of elements of the set $\widehat{\mathbb{N}}$ cannot be
expressed within a numeral system using only \G1. It is necessary
to introduce in a reasonable way a more powerful numeral system
and to define   new numerals (for instance, \ding{173},
\ding{174}, etc.) of this system that would allow one to fix the
set (or sets) somehow. In general, due to Postulate~1 and~2, for
any fixed numeral  system $\mathcal{A}$ there always be sets that
cannot be described using $\mathcal{A}$.

Let us give one more example illustrating properties of grossone.

\begin{example}
\label{e3_Lagrange}

Analogously to (\ref{4.1}), the set, $\mathbb{E}$, of even natural
numbers can be written now in the form
 \beq
\mathbb{E} = \{ 2,4,6 \hspace{5mm} \ldots  \hspace{5mm} \G1-4,
\hspace{2mm}\G1-2, \hspace{2mm} \G1 \}. \label{4.1.0}
       \eeq
Due to Postulate 3 and the IUA (see \cite{Sergeyev,informatica}),
it follows that the number of elements of the set of even numbers
is equal to $\frac{\G1}{2}$ and \G1 is even. Note that the next
even number is $\G1+2$ but it is not natural. In fact, since
$\G1+2  > \G1$, it is extended natural (see (\ref{4.2.2})). Thus,
we can write down not only initial (as it is done traditionally)
but also the final part of (\ref{4.4.1})
  \[
\begin{array}{cccccccccc}
 2, & 4, & 6, & 8,  & 10, & 12, & \ldots  &
\G1 -4,  &    \G1  -2,   &    \G1    \\
 \updownarrow &  \updownarrow & \updownarrow  &
\updownarrow  & \updownarrow  &  \updownarrow  & &
  \updownarrow    & \updownarrow   &
  \updownarrow
   \\
 1, &  2, & 3, & 4 & 5, & 6,   &   \ldots  &    \frac{\G1}{2} - 2,  &
     \frac{\G1}{2} - 1,  &    \frac{\G1}{2}   \\
     \end{array}
\]
concluding so (\ref{4.4.1})   in a complete accordance with
Postulate~3.

Suppose now that we have a set $A$ that has $k$ elements and all
its elements are multiplied by a constant in order to form the set
$B$. Then the number of the elements of the resulting set $B$ will
be the same as in the initial set $A$ independently on the fact
whether $k$ is finite or infinite. For instance, if we take
$A=\mathbb{N}$ then it has grossone elements. By choosing the set
$B=\{y: y=2x, x \in \mathbb{N} \}$, we have (see (\ref{4.1})) that
\[
B=\{ 2, 4, 6,   8,      \ldots \G1-4, \G1  -2, \G1, \G1  +2,
\G1+4, \ldots 2\G1-4, 2\G1  -2, 2\G1 \},
\]
i.e., it also has grossone elements. All elements of the set $B$
are even. Numbers $2, 4, 6,   8,      \ldots \G1-4, \G1  -2, \G1$
are even natural numbers and $\G1 +2, \G1+4, \ldots 2\G1-4, 2\G1
-2, 2\G1$ are even extended natural numbers.
 \hfill $\Box$
\end{example}

It is worth  noticing that the new numeral system allows us to
avoid many other `paradoxes' related to infinities and
infinitesimals (see \cite{Sergeyev,informatica,Korea}). For
instance, let us return to Hilbert's paradox of the Grand Hotel
presented in Section~\ref{s1}. In the original formulation of the
paradox, the number of   rooms in the Hotel is countable. In our
terminology, such a definition is not sufficiently precise. It is
necessary to indicate explicitly the infinite number of rooms in
the Hotel. Suppose that it has \G1 rooms. When a new guest
arrives, it is proposed to move the guest occupying room 1 to room
2, the guest occupying room 2 to room 3, etc. Finally, the guest
from room \G1 should be moved to room \G1+1 but the Hotel has only
\G1 rooms. As a result, the person from the last room should leave
the Hotel.

Thus, when the Hotel is full, no more new guests can be
accommodated in it if one wants that all guests living in   the
Hotel before the arrival of the newcomer remain inside. This
result corresponds perfectly to Postulate~3 and to the situation
taking place in hotels with a finite number of rooms.

Let us consider now the issue regarding a more systematic way to
produce numerals including \G1. In order to express more numbers
having finite, infinite, and infinitesimal parts, records similar
to traditional positional numeral systems can be used (see
\cite{Sergeyev,informatica}). To construct a number $C$ in the new
numeral positional system with the base~\G1, we subdivide $C$ into
groups corresponding to powers of \G1:
 \beq
  C = c_{p_{m}}
\G1^{p_{m}} +  \ldots + c_{p_{1}} \G1^{p_{1}} +c_{p_{0}}
\G1^{p_{0}} + c_{p_{-1}} \G1^{p_{-1}} + \ldots   + c_{p_{-k}}
 \G1^{p_{-k}}.
\label{3.12}
       \eeq
 Then, the numeral
 \beq
  C = c_{p_{m}}
\G1^{p_{m}}    \ldots   c_{p_{1}} \G1^{p_{1}} c_{p_{0}}
\G1^{p_{0}} c_{p_{-1}} \G1^{p_{-1}} \ldots c_{p_{-k}}
 \G1^{p_{-k}}
 \label{3.13}
       \eeq
represents  the number $C$, where all numerals $c_i$ are expressed
in a traditional numeral system we are used to express finite
numbers and are called \textit{grossdigits}. They express finite
positive or negative numbers (i.e., all $c_i\neq0$) and show how
many corresponding units $\G1^{p_{i}}$ should be added or
subtracted in order to form the number $C$.

Numbers $p_i$ in (\ref{3.13}) are  sorted in the decreasing order
with $ p_0=0$
\[
p_{m} >  p_{m-1}  > \ldots    > p_{1} > p_0 > p_{-1}  > \ldots
p_{-(k-1)}  >   p_{-k}.
 \]
They are called \textit{grosspowers} and they themselves can be
written in the form (\ref{3.13}).
 In the record (\ref{3.13}), we write
$\G1^{p_{i}}$ explicitly because in the new numeral positional
system  the number   $i$ in general is not equal to the grosspower
$p_{i}$. This gives the possibility to write down numerals without
indicating grossdigits equal to zero.

The term having $p_0=0$ represents the finite part of $C$ because,
due to (\ref{3.2.1}), we have $c_0 \G1^0=c_0$. The terms having
finite positive gross\-powers represent the simplest infinite
parts of $C$. Analogously, terms   having   negative finite
grosspowers represent the simplest infinitesimal parts of $C$. For
instance, the  number $\G1^{-1}=\frac{1}{\G1}$ is infinitesimal.
It is the inverse element with respect to multiplication for \G1:
 \beq
\G1^{-1}\cdot\G1=\G1\cdot\G1^{-1}=1.
 \label{3.15.1}
       \eeq
Note that all infinitesimals are not equal to zero. Particularly,
$\frac{1}{\G1}>0$ because it is a result of division of two
positive numbers. All of the numbers introduced above can be
grosspowers, as well, giving thus a possibility to have various
combinations of quantities and to construct  terms having a more
complex structure.

\begin{example}
\label{e4_Lagrange} In this example, it is shown   how to write
down numerals in the new positional numeral system and   how the
value of the number is calculated:
\[
 C_1=17.21\G1^{52.4\mbox{\small{\G1}}-72.1}\,
 134\G1^{81.43}
 7.02\G1^{0}52.1\G1^{-9.2}({\mbox{\small{-}}}0.23)\G1^{-3.7\mbox{\tiny\G1}}
=
 \]
 \[
17.21\G1^{52.4\mbox{\small{\G1}}-72.1}\,+
134\G1^{81.43}+7.02\G1^{0}+52.1\G1^{-9.2}-0.23\G1^{-3.7\mbox{\tiny\G1}}.
\]
The number $C_1$ above has two infinite parts of the type
$\G1^{52.4\mbox{\small{\G1}}-72.1}$ and $\G1^{81.43}$,   a finite
part corresponding to $\G1^{0}$, and two infinitesimal parts of
the type $\G1^{-9.2}$ and $\G1^{-3.7\mbox{\tiny\G1}}$. The
corresponding grossdigits show how many units of each kind should
be taken (added or subtracted) to form~$C_1$. \hfill$\Box$
\end{example}

\section{Numerical computations and modelling using the new methodology}

Let us start by considering what do we have instead of series when
we apply the new methodology, in particular, what happens in the
case of divergent series with alternating signs. As was already
mentioned,  the numeral $\infty$   is excluded from our new
mathematical language since we are able now to express explicitly
different infinite numbers. In fact, records of the type
$\sum_{i=1}^{\infty}a_i$ become  a kind of $\sum_{i=1}^{many}a_i$
and are not sufficiently precise. In order to define a sum
(independently on the fact whether the number of items in it is
finite or infinite), it is necessary  to indicate  explicitly how
many items we want to sum up. If the number of items in a sum is
infinite then, as it happens for the finite case, different
numbers of items in a sum lead  to different answers (that can be
infinite, finite, or infinitesimal). Let us give just two examples
(see \cite{informatica,Dif_Calculus} for a more detailed
discussion).

\begin{example}
\label{e5_Lagrange}

We start from the famous series
\[
S_1 = 1-1+1-1+1-1+\ldots
\]
 In literature, there exist  many
approaches giving different answers regarding the value of this
series (see \cite{Knopp}). All of them use various notions of
average to calculate the series. However, the notions of the sum
and of an average are two different things. In our approach, we do
not use the notion of series and do not appeal to an average. We
indicate explicitly the number of items, $k$, in the sum (where
$k$ can be finite or infinite) and calculate it directly:
\[
S_1(k)=\underbrace{1-1+1-1+1-1+1-\ldots}_{k} = \left \{
\begin{array}{ll} 0, &
  \mbox{if  } k=2n,\\
1, &    \mbox{if  } k=2n+1,\\
 \end{array} \right.
\]
and it is not important is  the number $k$ finite or infinite. For
example, for $k=2\G1$ we have $S_1(2\G1)=0$ and for $k=2\G1-1$ we
obtain $S_1(\mbox{2\G1}-1)=1$. \hfill$\Box$
\end{example}

It is important to emphasize that, as it happens in the case of
the finite number of items in a sum, the obtained answers do not
depend on the way the items in the entire sum are re-arranged. In
fact, if we  know the exact infinite number of items in the sum
and the order of alternating the signs is clearly defined, we know
also the exact number of positive and negative items in the sum.

Let us illustrate this point by supposing, for instance, that we
want to re-arrange the items in the sum $S_1(2\G1)$ in the
following way
\[
S_1(2\G1)=1+1-1+1+1-1+1+1-1+\ldots
\]
However, we know that the sum has $2\G1$ items and the number
$2\G1$ is even. This means that in the sum there are $\G1$
positive and $\G1$ negative items. As a result, the re-arrangement
considered above can continue only until the positive items
present in the sum will not finish and then it will be necessary
to continue to add only negative numbers. More precisely, we have
\[
S_1(2\G1)=\underbrace{1+1-1+1+1-1+\ldots+1+1-1}_{\G1 \mbox{
positive and }\frac{\mbox{\tiny{\G1}}}{2}\mbox{ negative
items}}\hspace{1mm}\underbrace{-1-1-\ldots-1-1-1}_{\frac{\mbox{\tiny{\G1}}}{2}\mbox{
negative items}}=0,
\]
where the result of the first part in this re-arrangement is
calculated as
$(1+1-1)\cdot\frac{\mbox{\tiny{\G1}}}{2}=\frac{\mbox{\tiny{\G1}}}{2}$
and the result of the second part is equal to
$-\frac{\mbox{\tiny{\G1}}}{2}$.

\begin{example}
\label{e6_Lagrange} Let us consider now the following divergent
series
 \[
S_2=  1 -2   + 3 -4 + \ldots
  \]
Again we should fix the number of items, $k$, in the sum $S_2(k)$.
Suppose that it contains grossone items. Then it follows
 \[
S_2(\G1) = 1 -2   + 3 -4 + \ldots - (\G1-2) + (\G1-1) - \G1 =
  \]
  \[
(\underbrace{1     + 3   + 5 + \ldots + (\G1-3) +
(\G1-1)}_{\frac{\tiny{\G1}}{2}\mbox{ items}} ) - ( \underbrace{2 +
4 + 6 + \ldots + (\G1-2) + \G1}_{\frac{\tiny{\G1}}{2}\mbox{
items}}  ) =
  \]
    \beq
   \frac{(1+\G1-1) \G1  }{4} -
\frac{(2+\G1) \G1  }{4} = \frac{ \G1^2 - 2 \G1  - \G1^2}{4}
 = - \frac{\G1}{2}.
 \label{Riemann15}
\eeq Obviously, if we change the number of items, $k$, then, as it
happens in the finite case, the results of summation will also
change. For instance, it follows $S_2(\G1-1)=\frac{\G1}{2}$ and
$S_2(\G1+1)=\frac{\G1}{2}+1$. \hfill$\Box$
\end{example}

Analogously to the passage from series to sums considered above,
we are able now to move from limits of expressions to the exact
evaluation of these expressions at points (finite, infinite or
infinitesimal) of our interest. Moreover, we can  calculate an
expression $f(x)$ independently on the fact of the existence of
the limit.   We are able  to change our way of thinking in sense
that instead of formulating problems in terms of limits by asking
`What does it happen when $x$ \textit{tends} to $\infty$?' we can
ask `What does it happen at different points of infinity?'

 Thus, limits are substituted by computation, at different points
 $x$, of
precise results $f(x)$ that can assume infinite, finite or
infinitesimal values and can be evaluated also in the cases where
limits do not exist. As a rule, the calculated values are
different for different infinite, finite, or infinitesimal values
of $x$. Note that the possibility of the direct evaluation of
expressions is very important (in particular, for automatic
computations) because it eliminates indeterminate forms from the
practice of computations.

For instance, in the traditional language if for a finite $a$,
$\lim_{x \rightarrow a}f(x)=0$ and $\lim_{y \rightarrow
\infty}g(y)=\infty$ then $\lim_{x \rightarrow a}f(x) \cdot \lim_{y
\rightarrow \infty}g(y)$ is an indeterminate form. In the new
language, this means that for any $x=a+z$ where $z$ is
infinitesimal, the value $f(a+z)$ is also infinitesimal and for
any infinite $y$ it follows that $g(y)$ is also infinite. In order
to be able to execute computations, we should behave ourselves as
we are used to do in the finite case. Namely, it is necessary to
choose $z$ and $y$, to evaluate $f(a+z)$ and $g(y)$. After
 we have performed  these operations it becomes possible to execute
multiplication $f(a+z)\cdot g(y)$ and to obtain the corresponding
result that can be infinite, finite or infinitesimal in dependence
of the values of $z$ and $y$ and the form of expressions $f(x)$
and $g(y)$.

 It
is possible also to execute other operations with infinitesimals
and infinities making questions with respect to $f(a+z)$ and
$g(y)$ that could not even  be  formulated using the traditional
language using limits. For instance, we can ask about the result
of the following expression
 \beq f(a+z_2)\left(\frac{g(y_1)}{f(a+z_1)}
-1.25 g(y_2)^3\right)
 \label{Lagrange_1}
 \eeq
 for two different infinitesimals $z_1, z_2$ and two different infinite
 values $y_1, y_2$.

\begin{example}
\label{e7_Lagrange} Let us consider an illustration regarding
computation of the product $f(a+z)\cdot g(y)$. For the sake of
simplicity we take $a=0$, $g(y)=y$, and
 \[
f(x) = \left\{ \begin{array}{ll} 2x, &   x  < 0,\\
1,& x=0,\\
x^3, & x>0. \end{array} \right.
 \]
If we want to calculate the product at points $z=\G1^{-1}$ and
$y=\G1$ then it follows
\[
f(a+z)\cdot g(y) = f(\G1^{-1})\cdot g(\G1)=\G1^{-3}\cdot \G1 =
\G1^{-2}.
\]
Analogously, $z=\G1^{-1}$ and $y=\G1^4$ give
\[
f(\G1^{-1})\cdot g(\G1^4)=\G1^{-3}\cdot \G1^4 = \G1^{1}
\]
and for $z=-2\G1^{-1}$ and $y=\G1$ we obtain
\[
f(-2\G1^{-1})\cdot g(\G1)=-4\G1^{-1}\cdot \G1 = -4.
\]
We end this example by calculating the result of the expression
(\ref{Lagrange_1}) for  $z_1=-2\G1^{-1}$, $z_2=-5\G1^{-4}$,
$y_1=\G1^2$, and $y_2=\G1$
\[
f(a+z_2)\left(\frac{g(y_1)}{f(a+z_1)}  -1.25 g(y_2)^3\right) =
f(-5\G1^{-4})\left(\frac{g(\G1^2)}{f(-2\G1^{-1})}-1.25
g(\G1)^3\right) =
\]
\[
\hspace{2mm}-10\G1^{-4}\cdot \left(\frac{\G1^2}{-4\G1^{-1}}
-1.25\G1^3\right)=-10\G1^{-4}\cdot\left(-0.25\G1^{3}-1.25\G1^3\right)=15\G1^{-1}.\hspace{2mm}\Box
\]
\end{example}

We conclude the paper by showing how the distinction between
mathematical objects and tools of their observation  helps us in
solving probabilistic questions and introduces the ideas of
relativity in Mathematics. In particular, we intend to show that
the new approach allows us to distinguish the impossible event
having the probability equal to zero (i.e., $P(\varnothing)=0$)
from those events that from the traditional point of view have the
probability equal to zero but can occur.

\begin{figure}[t]
  \begin{center}
    \epsfig{ figure = 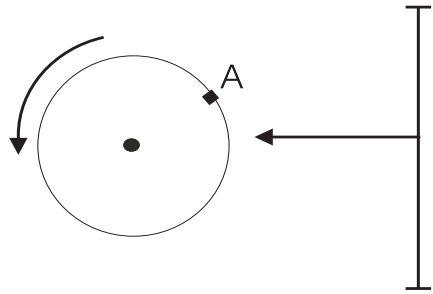, width = 2.2in, height = 1.5in,  silent = yes }
    \caption{What is the probability that the rotating disk stops in
such a way that the point $A$ will be exactly in front of the
arrow? }
 \label{Big_paper2}
  \end{center}
\end{figure}

 Let us consider the problem presented in
Fig.~\ref{Big_paper2} from the point of view of the traditional
probability theory. We start to rotate a disk having radius $r$
with the point $A$ marked at its border and we would like to know
the probability $P(E)$ of the following event $E$:  the disk stops
in such a way that the point $A$ will be exactly in front of the
arrow fixed at the wall.   Since the point $A$  is an entity that
has no extent, it is calculated by considering the following limit
\[
P(E) = \lim_{h \rightarrow 0}\frac{h}{2\pi r}=0.
\]
where $h$ is an arc of the circumference containing $A$ and $2\pi
r$ is its length.

However,   the point $A$ can stop in front of the arrow, i.e.,
this event is not impossible and its probability should be
strictly greater than zero, i.e., $P(E)>0$.  Obviously, this
example is a particular manifestation of the general fact that, if
$\xi$ is any continuous random value and $a$ is any real number
then $P(\xi=a)=0$.   While for a discrete random variable one
could say that an event with probability zero is impossible, this
can not be said in the terms of the traditional probability theory
for any continuous random variable.

Let us see what we can say with respect to this problem by using
the new methodology. The problem under consideration deals with
points located on the circumference $C$ of the disk. Thus, we need
a definition of the term `point'  and mathematical tools allowing
us to indicate a point on the circumference.  If we accept (as is
usually done in modern Mathematics) that  a \textit{point} is
determined by a numeral $x$ called the \textit{coordinate of the
point} where $x \in \mathcal{S}$ and
 $\mathcal{S}$ is a set of numerals,    then we can
indicate  the point  by its coordinate  $x$  and are able to
execute required calculations. The choice of the numeral system
$\mathcal{S}$ defines both the kind of numerals expressible in
this system and  the quantity (finite or infinite) of these
numerals (see \cite{Dif_Calculus,Korea} for a detailed
discussion). As a consequence, we are not able to work with those
points which coordinates are not expressible in the chosen numeral
system $\mathcal{S}$ (recall Postulate~2).

Different numeral systems can be chosen to express coordinates of
the points in dependence  on the precision level we want to
obtain. In some sense, the situation with counting points is
similar to the work with a microscope:  we decide the level of the
precision we need and obtain a result   dependent on the chosen
level. If we need a more precise or a more rough answer, we change
the level of the accuracy of our microscope. In the moment when we
have have decided which lens (numeral system) we put in the
microscope we decide which objects (points, arcs, etc.) we are
able to observe, to measure, and to work with.

The formalization of the concept `point' introduced above allows
us to execute  more accurate  computations having, as it always
happens in any process of the measurement, their own accuracy.
Suppose that we have chosen a numeral system $\mathcal{S}$
allowing one to observe $K$  points on the circumference.
Definition of the notion \textit{point} allows us to define
elementary events in our experiment as follows:  the disk has
stopped and the arrow indicates a point. As a consequence, we
obtain that the number, $N(\Omega)$, of all possible elementary
events, $e_i$,  in our experiment is equal to $K$ where
$\Omega=\cup_{i=1}^{N(\Omega)}e_i$ is the sample space of our
experiment. If our disk is well balanced, all elementary events
are equiprobable and, therefore, they have the same probability
equal to $\frac{1}{N(\Omega)}$  and the accuracy of any further
computation with this probabilistic model will be equal to
$\frac{1}{N(\Omega)}$. Thus,   we can calculate $P(E)$ directly by
subdividing the number, $N(E)$, of favorable elementary events by
the number, $K=N(\Omega)$, of all possible events.

For example, if we
 use numerals of the type
$\frac{i\cdot 2\pi r}{\G1}, i \in \mathbb{N},$ then $K=\G1$ and,
since the number of the points is infinite and the length of the
circumference is finite, our points are infinitesimally close,
i.e., the probabilistic model is continuous. The chosen numerals
define   the accuracy of the model and do not allow us to answer
to questions regarding objects having an extension on the
circumference that is less than $\frac{2\pi r}{\G1}$.

The number $N(E)$ depends on our decision about how many numerals
we want to use to represent the point $A$. If we decide that the
point~$A$ on the circumference is represented by $m$ numerals  we
obtain
\[
P(E) = \frac{N(E)}{N(\Omega)}= \frac{m}{K} = \frac{m}{\G1} > 0.
\]
where the number $\frac{m}{\G1}$ is infinitesimal if $m$ is
finite. Note that this representation is very interesting also
from the point of view of distinguishing   the notions `point' and
`arc'. When $m$ is finite than we deal with a point, when $m$ is
infinite we deal with an arc.

In the case we need the probabilistic model with a higher
accuracy, we can choose, for instance, numerals of the type $
\frac{i \cdot 2\pi r}{\G1^{-2}} , 1 \le i \le \G1^{2},$ for
expressing points on the circumference. In this way we also obtain
a continuous model with the order that is higher than in the
previous case. It follows $K=\G1^2$ and for a finite $m$ we obtain
the infinitesimal probability $P(E) = \frac{m}{\G1^2} > 0$.

In contrast, if we need a rough probabilistic model and decide to
work with   a finite number, $K$, of points on the circumference,
then we have the discrete model. In this case, the probability
$P(E)$ will be finite,  and the model does not allow us to answer
to questions regarding objects having an extension on the
circumference that is less than~$\frac{2\pi r}{K}$.

As we have shown by the example above, in our approach, for both
cases, the discrete and  the continuous one, only the impossible
event has the probability equal to zero. All other events have
positive probabilities that can be finite or infinitesimal in
dependence of the accuracy of the chosen probabilistic model.
Thus, the obtained probabilities are not absolute, i.e., there is
again a straight analogy with Physics where  results of the
observation   have a precision  determined by the used instrument.
Moreover, the new approach allows us to look at the same
mathematical object (like it happens in Physics for physical
objects) as continuous or discrete in dependence on the chosen
instrument of the observation (see \cite{Dif_Calculus} for a
detailed discussion related to this issue).

\bigskip


%


\normalsize


\footnotesize
\begin{thebibliography}{99}


%



\bibitem{Cantor} \textsc{Cantor~G.},  \textit{Contributions to the Founding of the
Theory of Transfinite Numbers}, Dover Publications, New York 1955.


\bibitem{Whorf} \textsc{Carroll~J.B. (Ed.)},
 \textit{Language, Thought, and Reality: Selected Writings of
  {B}enjamin {L}ee {W}horf}, MIT Press, 1956.


\bibitem{Cauchy}   \textsc{Cauchy~A.L.},
  \textit{Le Calcul infinit\'esimal}, Paris 1823.

\bibitem{Conway} \textsc{Conway~J.H. and Guy~R.K.},  \textit{The Book of
 Numbers},
Springer-Verlag, New York  1996.


\bibitem{DAlembert}  \textsc{d'Alembert~J.},
 \textit{Diff\'erentiel},
  Encyclop\'edie, ou dictionnaire raisonn\'e des sciences,
  des arts et des m\'etiers, \textbf{4}, Paris (1754).

\bibitem{Godel}  \textsc{G\"{o}del~K.},
    \textit{The Consistency of the {Continuum-Hypothesis}},
 Princeton University Press, Princeton 1940.

\bibitem{Godel_1931}  \textsc{G\"{o}del~K.},
  \textit{\"{U}ber formal unentscheidbare {S}\"{a}tze der
  {P}rincipia {M}athematica und verwandter {S}ysteme},
 Monatshefte f\"{u}r Mathematik und Physik \textbf{38} (1931),
 173--198.

\bibitem{Gordon} \textsc{Gordon~P.},  \textit{Numerical Cognition without Words:
 {E}vidence from {A}mazonia}, Science  \textbf{306} (15 October)  (2004),
 496--499.

\bibitem{Hardy}  \textsc{Hardy~G.H.}, \textit{Orders of infinity},
  Cambridge  University Press, Cambridge 1910.


\bibitem{Hilbert}  \textsc{Hilbert~D.},
 \textit{Mathematical Problems: Lecture delivered before the {I}nternational
  {C}ongress of {M}athematicians at {P}aris in 1900}, Bulletin of
   the American Mathematical Society \textbf{8} (1902), 437--479.

\bibitem{Knopp} \textsc{Knopp~K.},  \textit{Theory and Application
of Infinite Series}, Dover Publications, New York 1990.

\bibitem{Leibniz} \textsc{Leibniz~G.W.  and  Child~J.M.},  \textit{The Early
Mathematical Manuscripts of {L}eibniz},    Dover Publications, New
York 2005.


\bibitem{Mayberry} \textsc{Mayberry~J.P.},  \textit{The Foundations of Mathematics
 in the Theory of Sets}, Cambridge Univ. Press, Cambridge  2001.


\bibitem{Newton} \textsc{Newton~I.},  \textit{Method of
Fluxions},  1671.



\bibitem{Pica} \textsc{Pica~P.,   Lemer~C.,  Izard~V.,  Dehaene~S.},
 \textit{Exact and Approximate Arithmetic in an Amazonian Indigene Group},
 Science \textbf{306} (15 October) (2004), 499--503.

\bibitem{Robinson} \textsc{Robinson~A.},  \textit{Non-Standard Analysis},
 Princeton Univ. Press, Princeton  1996.

\bibitem{Sapir}  \textsc{Sapir~E.},
 \textit{Selected Writings of {E}dward {S}apir in Language,
  Culture and Personality}, University of California Press,
 Princeton 1958.

\bibitem{Sergeyev} \textsc{Sergeyev~Ya.D.},  \textit{Arithmetic of
Infinity}, Edizioni Orizzonti Meridionali, CS 2003.

\bibitem{chaos} \textsc{Sergeyev~Ya.D.},  \textit{Blinking Fractals and Their
 Quantitative Analysis Using Infinite and Infinitesimal Numbers},
 Chaos, Solitons  $\&$  Fractals. \textbf{33} 1 (2007), 50--75.


\bibitem{informatica} \textsc{Sergeyev~Ya.D.}, \textit{A New Applied Approach
for Executing Computations with Infinite and Infinitesimal
 Quantities}, Informatica  \textbf{19} 4 (2008), 567--596.

\bibitem{Menger} \textsc{Sergeyev~Ya.D.}, \textit{Evaluating the Exact
Infinitesimal Values of Area of Sierpinski's Carpet and Volume of
{M}enger's Sponge}, Chaos, Solitons  $\&$  Fractals  \textbf{42} 5
(2009), 3042--3046.


\bibitem{Dif_Calculus} \textsc{Sergeyev~Ya.D.}, \textit{Numerical  Point of View on Calculus for Functions Assuming Finite,
 Infinite, and Infinitesimal Values Over Finite, Infinite,
 and Infinitesimal Domains}, Nonlinear Analysis Series A:
 Theory, Methods $\&$ Applications  \textbf{71} 12 (2009), e1688--e1707

\bibitem{Korea} \textsc{Sergeyev~Ya.D.}, \textit{Numerical Computations
and Mathematical Modelling with Infinite and Infinitesimal
Numbers}, J. Applied Mathematics $\&$ Computing  \textbf{29}
(2009), 177--195.

\bibitem{Sergeyev_patent} \textsc{Sergeyev~Ya.D.},
\textit{Computer System  for Storing  Infinite, Infinitesimal, and Finite
Quantities and Executing Arithmetical
 Operations with Them}, EU patent 1728149 (2009).



\bibitem{first} \textsc{Sergeyev~Ya.D.}, \textit{Counting Systems and
the {F}irst {H}ilbert Problem}, Nonlinear Analysis Series A:
Theory, Methods $\&$ Applications  \textbf{72}  3-4 (2010),
1701--1708.



\bibitem{Sergeyev_Garro} \textsc{Sergeyev~Ya.D. and Garro~A.},   \textit{Observability of
Turing Machines:  A Refinement of the Theory of Computation}
 Informatica (2010) (in press).




\bibitem{www} \textit{The Infinity Computer web page}, \url{http://www.theinfinitycomputer.com}




\end{thebibliography}
\end{document}